\documentclass[11pt]{amsart}
\usepackage{amsmath, amsthm, latexsym, amssymb, amsfonts, epsfig, epsf}

\newenvironment{pf}{\proof[\proofname]}{\endproof}
\theoremstyle{plain}
\newtheorem{Th}{Theorem}[section]
\newtheorem{Cor}[Th]{Corollary}
\newtheorem{Conj}[Th]{Conjecture}
\newtheorem{Prop}[Th]{Proposition}
\newtheorem{Lemma}[Th]{Lemma}
\numberwithin{equation}{section} \theoremstyle{definition}
\newtheorem{Rem}[Th]{Remark}

\newtheorem{Ex}[Th]{Example}

\newtheorem{Def}[Th]{Definition}

\newcommand{\C}{\mathbb C}
\newcommand{\mS}{\mathbb S}

\newcommand{\R}{\mathbb R}

\newcommand{\A}{\mathbb A}
\newcommand{\G}{\Gamma}

\newcommand{\D}{\Delta}

\newcommand{\conv}{\operatorname{conv}}

\newcommand{\rs}[1]{Section~\ref{S:#1}}
\newcommand{\rl}[1]{Lemma~\ref{L:#1}}
\newcommand{\rp}[1]{Proposition~\ref{P:#1}}

\newcommand{\rex}[1]{Example~\ref{ex:#1}}
\newcommand{\re}[1]{(\ref{e:#1})}
\newcommand{\rcj}[1]{Conjecture \ref{Cj:#1}}
\newcommand{\rc}[1]{Corollary~\ref{C:#1}}
\newcommand{\rt}[1] {Theorem~\ref{T:#1}}

\newcommand{\rf}[1]{Figure~\ref{F:#1}}

\begin{document}


\title{Bezout Inequality for Mixed volumes}
\author{Ivan Soprunov}
\address[Ivan Soprunov]{Department of Mathematics\\ Cleveland State University\\ Cleveland, OH USA}
\email{i.soprunov@csuohio.edu}
\author{Artem Zvavitch}
\address[Artem Zvavitch]{Department of Mathematical Sciences\\ Kent State University\\ Kent, OH USA}
\email{zvavitch@math.kent.edu}
\thanks{The first author is  supported in part by NSA Grant H98230-13-1-0279. The second author is supported in part by U.S. National Science Foundation Grant DMS-1101636  and by the  Simons Foundation.}

\begin{abstract} In this paper we consider the following  analog of Bezout inequality for mixed volumes:
$$V(P_1,\dots,P_r,\D^{n-r})V_n(\D)^{r-1}\leq \prod_{i=1}^r V(P_i,\D^{n-1})\ \text{ for }2\leq r\leq n.$$
We show that the above inequality is true when $\D$ is an $n$-dimensional simplex and $P_1, \dots, P_r$ are convex bodies in $\R^n$.  We conjecture that if the above inequality 
is true for all convex bodies $P_1, \dots, P_r$, then $\D$ must be an $n$-dimensional simplex.  We prove that if the above inequality is true  for all convex bodies $P_1, \dots, P_r$, then $\D$ must be indecomposable (i.e. cannot be written as the Minkowski sum of two convex bodies which are not homothetic to $\D$), which confirms the conjecture 
when $\D$ is a simple polytope and in the $2$-dimensional case. Finally, we connect the inequality to an inequality on the volume of orthogonal projections of convex bodies as well as prove an isomorphic version of the inequality.

\end{abstract}

\maketitle


\section{Introduction}
The classical Bezout inequality in algebraic geometry relates the degrees of 
hypersurfaces to the degree of their intersection. More precisely, let $X$ be a closed
algebraic set in an affine space $\A^n$ over the field of complex numbers $\C$
(or an algebraically closed  field). Its degree
$\deg X$ is defined as the number of intersection points of $X$ with a generic affine subspace
of complementary dimension. Given $1\leq r\leq n$ hypersurfaces $X_1,\dots, X_r$ in $\A^n$
whose intersection has pure codimension $r$, the Bezout
inequality says 
$$\deg(X_1\cap\dots\cap X_r)\leq \deg X_1\cdots \deg X_r.$$
(see, for instance, Proposition 8.4 of \cite{F} and examples therein).
The theory of Newton polytopes provides a beautiful interconnection
 between algebraic geometry and convex geometry. One of its central results is the
Bernstein--Kushnirenko--Khovanskii theorem which expresses the number of intersection points
of $n$ hypersurfaces in $(\C\setminus\{0\})^n$ with fixed Newton polytopes and generic coefficients as the
normalized mixed volume of the Newton polytopes, see \cite{Be, Kush, Kho}. In particular, if $X$ is a generic hypersurface
with Newton polytope $P$ and $H$ is a line given by $n-1$ generic linear forms, the degree of $X$ equals
$$\deg X =\#(X\cap H)=n!V(P,\D_n,\dots,\D_n),$$
where $\D_n=\conv\{0,e_1,\dots,e_n\}$ is the standard $n$-simplex (the Newton polytope of a generic linear form),
and $V(A_1,\dots,A_n)$ denotes the mixed volume of convex bodies $A_1,\dots,A_n$.

Similarly, 
the degree of the intersection $X_1\cap\dots\cap X_r$ is given by
$$\deg (X_1\cap\dots\cap X_r) =n!V(P_1,\dots,P_r,\D_n,\dots,\D_n),$$
where $P_i$ is the Newton polytope of $X_i$. In the mixed volume above, the simplex $\D_n$ 
is repeated $n-r$ times. In what follows we abbreviate this as $V(P_1,\dots,P_r,\D_n^{n-r})$,
using standard notation from the theory of mixed volumes.

Therefore, the Bezout inequality becomes the inequality for mixed volumes:
$$n!V(P_1,\dots,P_r,\D_n^{n-r})\leq \prod_{i=1}^r n!V(P_i,\D_n^{n-1}).$$
As the volume of $\D_n$ is $1/n!$ we can rewrite this as
\begin{equation}\label{e:Bl}
V(P_1,\dots,P_r,\D_n^{n-r})V_n(\D_n)^{r-1}\leq  \prod_{i=1}^r V(P_i,\D_n^{n-1}).
\end{equation}

We call this inequality {\it the Bezout Inequality for mixed volumes}, and the purpose of this paper is to study this inequality and its offspring and applications. 

In \rs{ineq} we start by
giving a geometric proof of the inequality which works for 
arbitrary convex bodies $P_1,\dots,P_r$ and arbitrary simplex $\D_n$. In addition, we
give a necessary and sufficient condition when the inequality becomes equality when
the $P_i$ are convex polytopes.
It is thus natural to ask if the class of simplices  
is a unique class of convex bodies for which the above inequality is true
when $r\geq 2$. More precisely, fix a convex body $D$ and assume
\begin{equation}\label{e:BE}
V(K_1,\dots, K_r,D^{n-r})V_n(D)^{r-1}\leq  \prod_{i=1}^r V(K_i,D^{n-1})
\end{equation}
holds for all  convex bodies $K_1,\dots, K_r$ in $\R^n$. What can be said about $D$?

In \rs{others} we prove that $D$ must be indecomposable. In particular, if $D$ is a
simple polytope then $D$ must be an $n$-simplex. Based on 
this observation, as well as other results in the paper, we propose the following conjecture.

 \begin{Conj}\label{Cj:BE} Fix  integers $2\leq r\leq n$ and 
 let $D$ be an $n$-dimensional convex body which satisfies the Bezout inequality \re{BE} for
 all convex bodies $K_1,\dots, K_r$ in $\R^n$. Then $D$ is an $n$-simplex.
 \end{Conj}
 
 The indecomposability result of \rt{decompose} gives us that the above conjecture is true in the case of $n=2$ (see \rs{others}). We also show in \rex{octahedron}, that there are indecomposable bodies $D$ in dimension $3$ and higher for which  \re{BE} is false. In \rs{proj} we connect inequality (\ref{e:BE}) to the inequalities related 
 to the volume of orthogonal projections of convex bodies. This connection helps us to provide more examples in support of Conjecture \ref{Cj:BE}, including that the body $D$ cannot be rotation invariant.
In \rs{iso} we provide an isomorphic version of (\ref{e:BE}), that is we show that there exists a constant $c_{n,r}$, depending on $n$ and $r$ only, such that 
\begin{equation}\label{e:isoBE}
V(K_1,\dots, K_r,D^{n-r})V_n(D)^{r-1}\leq c_{n,r}  \prod_{i=1}^r V(K_i,D^{n-1})
\end{equation}
is true for all convex bodies $K_1,\dots, K_r$ and $D$ in $\R^n$. Using the results of 
Fradelizi, Giannopoulos, Hartzoulaki, Meyer, and Paouris in \cite{FGM, GHP}
we show that $c_{n,r}=r^{r}/r!$ 
if $K_1,\dots, K_r$ are zonoids and 
$c_{n,r}\leq n^rr^{r}/r!$ in the case of general convex bodies. 
Finally, in \rs{dim2} we give a direct geometric proof
of inequality (\ref{e:isoBE}) in the
two dimensional case for the class of general convex bodies  
with the optimal constant $c_{2,2}=2$. It was pointed out to us by Christos Saroglou that
this 2-dimensional result was also proved by S.~Artstein-Avidan, D.~Florentin, and Y.~Ostrover
in \cite[Proposition 5.1]{AAFO}. 

\subsection*{Acknowledgment}  We are indebted to Christos Saroglou 
 for many valuable discussions and to the referee for helpful comments.

\section{Bezout Inequality for Mixed volumes}\label{S:ineq}

We start by giving basic definitions and setting up notation. As a general reference on the theory of 
convex sets and mixed volumes we use R. Schneider's book ``Convex bodies: the Brunn-Minkowski theory" \cite{Sch1}.

For $x,y\in\R^n$,  $x\cdot y$ denotes the inner product of $x$ and $y$. If $X$ and $Y$
are sets in $\R^n$ then $X+Y=\{x+y\ |\ x\in X, y\in Y\}$
is the {\it Minkowski sum} of $X$ and~$Y$. A {\it convex body} is a non-empty convex compact set. 
For a convex body $K$ the function $h_K(u)=\max\{x\cdot u\ |\ x \in K\}$ is the {\it support function} of $K$.
A {\it (convex) polytope} is the convex hull of a finite set of points. An $n$-dimensional polytope is called {\it simple} if
every its vertex is adjacent to exactly $n$ edges. A polytope which is the Minkowski sum 
of finitely many line segments is called a {\it zonotope}. Limits of zonotopes in the Hausdorff metric are called
{\it zonoids}, see \cite[Section 3.2]{Sch1} for details.
Let $V(K_1,\dots,K_n)$ denote the $n$-dimensional mixed volume of
$n$ convex bodies $K_1,\dots, K_n$ in $\R^n$.  We will also denote 
by $V_{n}(K)$ the $n$-dimensional (Euclidean) volume of $K$.

In this section we give a geometric proof of inequality (\ref{e:Bl}) and describe when the inequality becomes equality.
We need the following definition.

\begin{Def} A collection of $m$ convex sets in $\R^n$ is called {\it essential}
if the Minkowski sum of any $k$ of them has dimension at least $k$, for $1\leq k\leq m$.
\end{Def}

It is well known that  $n$ convex bodies in $\R^n$  form an essential collection if and only if their mixed volume is non-zero. This result goes back to Minkowski's work \cite{Min} (see also \cite[Theorem 5.1.7]{Sch1}).

 In the next statement and in what follows we use $A^u$ to denote the face of a convex body $A$
defined by the support hyperplane to $A$ with normal vector  $u$.

\begin{Prop}\label{P:mix}
Let  $A_1,\dots, A_n$ be convex polytopes in $\R^n$ and $K\subseteq A_1$, a convex body.
Then $V(K,A_2,\dots, A_n)=V(A_1,A_2,\dots, A_n)$ if and only if $K$ intersects every face $A_1^u$
for $u$ in the set
$$E=\{u\in \mS^{n-1}\ |\ (A_2^u,\dots, A_n^u)\ \text{\rm  is essential}\}.$$   
\end{Prop}

\begin{pf} We use the inductive formula for the mixed volume \cite[Theorem 5.1.6]{Sch1}. We have 
$$V(K,A_2,\dots, A_n)=\frac{1}{n}\int_{u\in \mS^{n-1}}h_K(u)\,dS(A_2,\dots, A_n,u),$$
where $h_K(u)$ is the support function of $K$ and $S(A_2,\dots, A_n,u)$ is the mixed area 
measure on the unit sphere $\mS^{n-1}$.

Note that the set $E$ is finite. Indeed, the Minkowski sum $A_2^u+\dots+ A_n^u$ is $(n-1)$-dimensional, hence, must
be a facet of $A'=A_2+\dots+A_n$. But each such facet corresponds to a unique $u\in \mS^{n-1}$. 
By the above, we can replace the right hand side with a finite sum
$$V(K,A_2,\dots, A_n)=\frac{1}{n}\sum_{u\in E}h_K(u)V(A_2^u,\dots, A_n^u).$$
Here $V(A_2^u,\dots, A_n^u)$ is the $(n-1)$-dimensional mixed volume of $A_2^u,\dots, A_n^u$
translated to the orthogonal subspace $u^\perp$.

Clearly, $h_K(u)\leq h_{A_1}(u)$ for any $u\in \mS^{n-1}$, 
as $K\subseteq A_1$.
Since $V(A_2^u,\dots, A_n^u)>0$ for $u\in E$,
 the equality $V(K,A_2,\dots, A_n)=V(A_1,A_2,\dots, A_n)$ holds if and only if
$h_K(u)=h_{A_1}(u)$ for all $u\in E$. The latter means that $K^u\cap A_1^u$ is non-empty.\end{pf}

\begin{Cor}\label{C:cor} 
Let $P$ be a convex polytope and $K\subseteq P$ a convex body.
Then $V(K,P^{n-1})=V_n(P)$ if and only if $K$ intersects every facet of $P$.
\end{Cor}

\begin{pf} Indeed, the set $E=\{u\in \mS^{n-1}\ |\ (P^u,\dots, P^u)\ \text{\rm  is essential}\}$
consists of exactly the normal vectors to the facets of $P$. 
\end{pf}

\begin{Rem}
It is a natural question to ask what necessary and sufficient conditions are required for 
$V(K,A_2,\dots, A_n)=V(A_1,A_2,\dots, A_n)$ to hold when $K\subset A_1$ and $A_2,\dots, A_n$ are
arbitrary convex bodies. This is an open question in general. See a detailed discussion in \cite[page 277]{Sch1}, 
as well as Conjecture 6.6.13, and Theorem~6.6.16 therein.
\end{Rem}

Let $K$ be a convex body and $\D$ an $n$-simplex in $\R^n$. We say that $K$ is in {\it convenient position
with respect to $\D$} if $K\subset\D$ and $K$ has a non-empty intersection with every facet of $\D$.
In general let $\lambda>0$ be
the largest number such that $\lambda K+v$ is contained in $\D$ for some $v\in\R^n$.
Then it is easy to see that 
 $\lambda K+v$ is in convenient position with respect to $\D$.

The next theorem is a generalization of \re{Bl} to arbitrary convex bodies and an arbitrary simplex.

\begin{Th}\label{T:ineq}
Let $K_1,\dots, K_r$ be convex bodies in $\R^n$, for $1\leq r\leq n$, and $\D$ an $n$-simplex.
Then the following inequality holds:
\begin{equation} \label{e:ineq}
V(K_1,\dots, K_r,\D^{n-r})V_n(\D)^{r-1}\leq  \prod_{i=1}^r V(K_i,\D^{n-1}).
\end{equation}
\end{Th}

\begin{pf} 
The inequality \re{ineq} is homogeneous in the $K_i$ and is independent 
of translations of the $K_i$ in $\R^n$. Therefore,  after a possible translation and dilation 
we may assume that all $K_i$ are contained in $\D$ and intersect each facet of $\D$, i.e.
each $K_i$ is in convenient position with respect to $\D$.
In this case, by \rc{cor}, we have  $V(K_i,\D^{n-1})=V_n(\D)$.
Now (\ref{e:ineq}) is equivalent to 
$$V(K_1,\dots, K_r,\D^{n-r})\leq V_n(\D),$$
which is true by the monotonicity of the mixed volume.
\end{pf}

Next we will give a description of  the $K_i$ for which the Bezout inequality  \re{ineq} becomes equality
in the case of polytopes. 

\begin{Th} Let $K_1,\dots, K_r$ for $2\leq r\leq n$ 
be convex polytopes in $\R^n$ which are in convenient position
with respect to a simplex~$\D$. Then equality in
the Bezout inequality \re{ineq} is attained at $K_1,\dots, K_r$ if and only if for any subset 
$\{K_{i_1},\dots, K_{i_s}\}\subseteq\{K_1,\dots, K_r\}$ the union $K_{i_1}\cup\dots\cup K_{i_s}$ has a non-empty intersection with every
$(n-s)$-dimensional face of~$\D$. 
\end{Th}

\begin{pf}
First, suppose we have equality in \re{ineq}. As in the proof of \rt{ineq} this is equivalent to
\begin{equation}\label{e:equiv}
V(K_1,\dots, K_r,\D^{n-r})= V_n(\D).
\end{equation}
By monotonicity of the mixed volume, it follows that for any subset $\{K_{i_1},\dots, K_{i_s}\}$
we also have
\begin{equation}\label{e:many}
V(K_{i_1},\dots, K_{i_s},\D^{n-s})= V_n(\D).
\end{equation}
We will show that in this case the union  $K_{i_1}\cup\dots\cup K_{i_s}$ intersects every
$(n-s)$-dimensional face of $\D$. Suppose not, so there exists a $(n-s)$-dimensional face $\G$ of 
$\D$ with $(K_{i_1}\cup\dots\cup K_{i_s})\cap\G=\varnothing$. Let $K$ be the convex hull of
$K_{i_1}\cup\dots\cup K_{i_s}$ which is contained in $\D$. 
Since $\G$ is a face of $\D$, it follows that $K\cap\G=\varnothing$, as well.
Choose $v\in \mS^{n-1}$ such that the hyperplane $H$ orthogonal to $v$ strictly separates $K$ and $\G$ and
$\G=\D^v$. Furthermore, let $\bar\D$ be the ``truncated" simplex, i.e. $\bar\D=\D\cap H^+$,
where $H^+$ is the half-space containing $K$. Then we have
\begin{equation}\label{e:monoton}
V(K_{i_1},\dots, K_{i_s},\D^{n-s})\leq V(\bar\D^s,\D^{n-s}),
\end{equation}
by monotonicity. Furthermore,  we have $h_{\bar\D}(u)\leq h_{\D}(u)$ for all $u\in \mS^{n-1}$
and $h_{\bar\D}(v) < h_{\D}(v)$ by construction. In addition, 
$$V\left((\bar\D^v)^{s-1}, (\D^v)^{n-s}\right)>0,$$
since $\dim\bar\D^v=n-1$ and $\dim\D^v=n-s$. Therefore,
$$V(\bar\D^s,\D^{n-s})=\frac{1}{n}\sum_{u\in \mS^{n-1}}h_{\bar\D}(u)V\left((\bar\D^u)^{s-1}, (\D^u)^{n-s}\right)<$$
$$\frac{1}{n}\sum_{u\in \mS^{n-1}}h_{\D}(u)V\left((\bar\D^u)^{s-1}, (\D^u)^{n-s}\right)=V(\bar\D^{s-1},\D^{n-s+1})\leq V_n(\D),$$
which together with \re{monoton} contradicts \re{many}.

For the other implication, suppose for any subset 
$\{K_{i_1},\dots, K_{i_s}\}\subseteq\{K_1,\dots, K_r\}$
the union $K_{i_1}\cup\dots\cup K_{i_s}$ has a non-empty intersection with every
$(n-s)$-dimensional face of $\D$. We claim that this implies \re{equiv}.
Clearly, this is true for $r=1$ since \re{ineq} is trivial for $r=1$.
We use induction on $r$.  By the inductive hypothesis,  the equalities \re{many} hold
whenever $s<r$. In particular, we have
\begin{equation}
V(K_1,\dots, K_{r-1},\D^{n-r+1})=V_n(\D).\nonumber
\end{equation}
Thus, it is enough to show that 
\begin{equation}\label{e:particular}
V(K_1,\dots, K_{r},\D^{n-r})=V(K_1,\dots, K_{r-1},\D^{n-r+1}),
\end{equation}
or, equivalently, 
\begin{equation}
\sum_{u\in \mS^{n-1}}h_{K_r}(u)V\!\left(K_1^u,\dots, K_{r-1}^u,(\D^u)^{n-r}\right)=
\sum_{u\in \mS^{n-1}}h_{\D}(u)V\!\left(K_1^u,\dots, K_{r-1}^u,(\D^u)^{n-r}\right).\nonumber
\end{equation}
Consider $u\in \mS^{n-1}$. If $K_r\cap\D^u\neq\varnothing$ then $h_{K_r}(u)=h_{\D}(u)$.
Otherwise, we claim that 
the mixed volume $V\!\left(K_1^u,\dots, K_{r-1}^u,(\D^u)^{n-r}\right)$ is zero and 
the above equality follows. Indeed, if
$\dim\D^u<n-r$ then the collection $\left(K_1^u,\dots, K_{r-1}^u,(\D^u)^{n-r}\right)$ is not essential. 
Assume $\dim\D^u=n-s$ for some $2\leq s\leq r$. By the condition of the theorem and since 
$K_r\cap\D^u=\varnothing$, the set $\{K_1,\dots,K_{r-1}\}$ can contain at most $s-2$ of the 
$K_i$ that are disjoint from $\D^u$. Therefore, there exist $r-s+1$ of the $K_i$ among 
 $\{K_1,\dots,K_{r-1}\}$ which intersect $\D^u$. We may assume they are
 $K_1,\dots, K_{r-s+1}$. But then the sum of $n-s+1$ bodies $K^u_1+\dots+K^u_{r-s+1}+(n-r)\D^u$
 has dimension $n-s$, which means that the collection $\left(K_1^u,\dots, K_{r-1}^u,(\D^u)^{n-r}\right)$ is not essential. 
 \end{pf}

\section{The case of other polytopes}\label{S:others}

In this section we study Conjecture \ref{Cj:BE}. We start with  two classical  statements describing the properties of log-concave sequences (see \cite{DP}, \cite{W} or \cite{H}).
\begin{Lemma}\label{L:1} 
Let $a_0\dots,a_n$ be a sequence of non-negative real numbers.
Then $$a_i^2\geq a_{i-1}a_{i+1}\text{ for all }1\leq i\leq n-1$$ if and only if 
$$a_ia_j\geq a_{i-1}a_{j+1}\text{ for all }1\leq i\leq j\leq n-1.$$
\end{Lemma}

\begin{pf} The ``if" statement is obvious. The ``only if" follows by induction on $j-i$.
\end{pf}

\begin{Lemma}\label{L:2} 
Let $a_0\dots,a_n$ be a sequence of non-negative real numbers.
Fix $0\leq m\leq n$ and let
$$C_i=\sum_{j=0}^m{m\choose j}a_{i+j},\quad 0\leq i\leq n-m.$$
Then 
\begin{equation}\label{e:L1}
a_i^2\geq a_{i-1}a_{i+1}\text{ for all }1\leq i\leq n-1
\end{equation} implies 
\begin{equation}\label{e:L2}
C_i^2\geq C_{i-1}C_{i+1}\text{ for all }1\leq i\leq n-m-1.
\end{equation} 
If, in addition, all the inequalities in \re{L2} are equalities then so are \re{L1}.
\end{Lemma}

\begin{pf} The proof is by induction on $m$. For $m=0$ the lines \re{L1} and \re{L2}
are the same. Let $C_i'=C_i+C_{i+1}$ for $0\leq i\leq n-m-1$. It is easy to see that 
$$C_i'=\sum_{j=0}^{m+1}{m+1\choose j}a_{i+j}.$$
For every $1\leq i\leq n-m-2$, by the inductive hypothesis and \rl{1}, we have
\begin{equation}\label{e:L3}
C_i^2\geq C_{i-1}C_{i+1},\quad C_iC_{i+1}\geq C_{i-1}C_{i+2},\quad C_{i+1}^2\geq C_{i}C_{i+2}.
\end{equation} 
This implies 
\begin{equation}\label{e:L4}
(C_i+C_{i+1})^2\geq (C_{i-1}+C_{i})(C_{i+1}+C_{i+2}),\text{ i.e. } C_i'^2\geq C_{i-1}'C_{i+1}'
\end{equation} 
and the first statement follows. If the inequality \re{L4} is, in fact, equality then
so are \re{L3} and we obtain the second statement of the lemma.
\end{pf}

We recall that the convex body $D \subset \R^n$ is called {\it indecomposable} if a representation 
$D=A+B$, where $A$ and $B$ are convex bodies, is only possible when $A$ and $B$ are homothetic to $D$ (see  \cite[Section 3.2]{Sch1}). 

\begin{Th}\label{T:decompose} Let $D$ be an $n$-dimensional convex body.
Suppose \re{BE} holds for all convex bodies $K_1,\dots, K_r$ in $\R^n$ where $2\leq r\leq n$. Then $D$ is indecomposable.
\end{Th}

\begin{pf} 
Suppose $D=A+B$ for some convex bodies $A$ and $B$.
We will show that $A$ is homothetic to $B$. 

First, we set $K_1=A$, $K_2=B$, and $K_i=D$ for $3\leq i\leq r$. Then the inequality \re{BE}
simplifies to 
\begin{equation}\label{e:simple}
V(A,B,D^{n-2})V_n(D)\leq  V(A,D^{n-1})V(B,D^{n-1}).
\end{equation}
Since $D=A+B$ we have 
$$V(A,B,D^{n-2})=\sum_{j=0}^{n-2}{n-2 \choose j}V(A^{n-1-j},B^{j+1}).$$
Similarly, 
$$V(A,D^{n-1})=\sum_{j=0}^{n-1}{n-1 \choose j}V(A^{n-j},B^{j}),\quad V(B,D^{n-1})=\sum_{j=0}^{n-1}{n-1 \choose i}V(A^{n-1-j},B^{j+1}),$$
and 
$$V_n(D)=\sum_{j=0}^{n}{n \choose j}V(A^{n-j},B^{j}).$$
Denoting $a_j=V(A^{n-j},B^j)$, the inequality \re{simple} becomes
$$
\left(\sum_{j=0}^{n-2}{n-2\choose j}a_{j+1}\right)\left(\sum_{j=0}^{n}{n \choose j}a_j\right)\leq 
\left(\sum_{j=0}^{n-1}{n-1 \choose j}a_j\right)\left(\sum_{j=0}^{n-1}{n-1 \choose j}a_{j+1}\right).
$$
Let 
$$C_i=\sum_{j=0}^{n-2}{n-2\choose j}a_{i+j},\quad 0\leq i\leq 2.$$
Then the above inequality states
\begin{equation}\label{e:BE2}
C_1(C_0+2C_1+C_2)\leq (C_0+C_1)(C_1+C_2).
\end{equation}
This is equivalent to $C_1^2\leq C_0C_2$.
Since the sequence $a_i$ satisfies the Alexandrov--Fenchel inequalities $a_i^2\geq a_{i-1}a_{i+1}$, by
\rl{2} we must have $C_1^2\geq C_0C_2$. Therefore, $C_1^2=C_0C_2$ which implies
 $a_i^2= a_{i-1}a_{i+1}$, again by \rl{2}.
 
 First, assume that $A$ and $B$ are $n$-dimensional.
 Notice that the equations $a_i^2= a_{i-1}a_{i+1}$ are homogeneous, so after possible rescaling we may
 assume that $B\subset A$. From the proof of \cite[Theorem 6.6.18]{Sch1} it follows that 
 $a_0=a_1$. Indeed, $a_1^2=a_0a_2$ is the case of equality in (6.6.52) from \cite[page 396]{Sch1}
 for $K=A$ and $L=B$ which gives $a_0=a_1$, see the end of the proof on \cite[page 372]{Sch1}.
  Together with $a_i^2= a_{i-1}a_{i+1}$ this implies that $a_0=\dots=a_n$.
 In particular, $a_0=a_n$, i.e. $A$ and $B$ have the same volume. Therefore, $A$ is homothetic to $B$.

 Now assume $\dim A=k< n$. Then $a_0=0$ and we obtain $a_0=\dots=a_n=0$. In particular, $\dim B< n$.
 Let $l=\dim B$. Since $D$ is $n$-dimensional we have  $k+l\geq \dim (A+B)=n$. But then the collection
 $(A^{n-l}, B^l)$ is essential (the Minkowski sum of any $m$ of them is at least $m$-dimensional for
  any $1\leq m\leq n$), which contradicts $a_l=V(A^{n-l},B^l)=0$.
\end{pf}

\begin{Cor} Suppose $D$ is a simple $n$-dimensional polytope satisfying  \re{BE}  for all bodies $K_1,\dots, K_r$ 
where $2\leq r\leq n$. Then $D$ is an $n$-simplex.
\end{Cor}

\begin{pf} This follows from the fact that any simple $n$-dimensional polytope is decomposable unless it is an $n$-simplex, see 
\cite[Theorem 15.1.4]{Gru2}.
\end{pf}

\begin{Cor} Suppose $D$ is a $2$-dimensional convex body in $\R^2$ satisfying 
$$V(K_1,K_2)V_2(D)\leq V(K_1,D)V(K_2,D)$$
for all convex bodies $K_1,K_2\subset\R^2$. Then $D$ is a $2$-simplex.
\end{Cor}
 
 \begin{pf} It is well known that any convex body $D$ in $\R^2$ is decomposable, unless $D$ is a
 segment or a $2$-simplex, see for example \cite[Theorem 3.2.11]{Sch1}.
\end{pf}
 
 The next example shows the existence of indecomposable polytopes for which \re{BE} fails.
 
 \begin{Ex}\label{ex:octahedron} Let $O=\{x\in\R^3\ |\ \sum_{i=1}^3|x_i|\leq 1\}$ be an octahedron in 
 $\R^3$. Then $O$ is indecomposable since all its 2-dimensional
 faces are triangles (see \cite[Corollary 3.2.13]{Sch1}). However, one can find $K_1$, $K_2$ such that \re{BE}
 fails with $r=2$ and $D=O$. Indeed, let $K_1$ be a segment connecting two opposite (non-intersecting) faces of $O$
 and $K_2$ one of those faces.  Let $O'$ (resp. $K_2'$) be the projection of $O$ (resp. $K_2$)
 onto the plane orthogonal to $K_1$. Then by (5.3.23)  of \cite[page 294]{Sch1} we obtain
 $$V(K_1,K_2,O)=\frac{V_1(K_1)}{3}V(K_2',O')\quad\text{and}\quad V(K_1,O,O)=\frac{V_1(K_1)}{3}V_2(O').$$
 Since $K_2'$ intersects all facets (sides) of $O'$ we have $V(K_2',O')=V_2(O')$, by \rc{cor}.
 On the other hand, $V(K_2,O,O)<V_3(O)$ since $K_2$ does not intersects all facets of $O$,
 again by \rc{cor}. Therefore, we obtain
 $$V(K_1,K_2,O)V_3(O)>V(K_1,O,O)V(K_2,O,O).$$

 \end{Ex}
 
 \section{Connections to Projections}\label{S:proj}

In this section we relate Conjecture \ref{Cj:BE} to inequalities involving orthogonal projections. 
In particular, we see how  \re{BE} may fail when we choose $K_1$ and $K_2$ to
be segments and $K_i=D$ for $3\leq i\leq r$. 

Let $K_1$ and $K_2$ be unit segments  and  $\xi, \nu \in \mS^{n-1}$ the corresponding direction vectors (chosen up to a sign). Then,  by (5.3.23)  of \cite[page 294]{Sch1},
$$
V(K_1, D^{n-1})=\frac{1}{n}V_{n-1}(D | \xi^\perp)\text{ and }V(K_2, D^{n-1})=\frac{1}{n}V_{n-1}(D | \nu^\perp),
$$
where $D | \xi^\perp$ denotes the orthogonal projection of  $D$ onto the hyperplane orthogonal to $\xi$. 
Let $U=K_1+K_2$. Similarly, for the orthogonal projection $D |  (\xi, \nu)^\perp$ of $D$ onto the $(n-2)$-dimensional subspace
$(\xi,\nu)^\perp$ we have 
$$
V_{n-2}(D |  (\xi, \nu)^\perp) V_2(U) =\binom{n}{2} V(U, U, D^{n-2}).
$$
By the linearity  of the mixed volume we get  $V(U, U, D^{n-2})=2 V(K_1, K_2, D^{n-2})$, and
$$
V_{n-2}(D |  (\xi, \nu)^\perp)V_2(U)  = n(n-1) V(K_1, K_2, D^{n-2}).
$$
Substituting the above calculations in \re{BE} we obtain
\begin{equation}\label{proj}
 \frac{n}{n-1}V_{n-2}(D |  (\xi, \nu)^\perp)V_2(U)V_n(D)\leq V_{n-1}(D | \xi^\perp) V_{n-1}(D | \nu^\perp).
\end{equation}
Equation (\ref{proj}) turns out to be quite useful to check the particular  cases of \rcj{BE}.  Direct calculations and the fact that the inequality \re{BE} is invariant under linear transformations shows that (\ref{proj}) (and thus \re{BE}) are false when $D$ is an $n$-dimensional ellipsoid  or parallelepiped (also both of those cases follow from \rt{decompose} as both are decomposable). 

We may also  generalize \rex{octahedron}, using (\ref{proj}). Indeed, consider the $n$-dimensional octahedron $O_n=\{x \in \R^n\ |\ \sum_{i=1}^n |x_i| \le 1\} $. Then $V_n(O_n)=2^n/n!$. Choose 
$$
\xi =({1}/{\sqrt{2}}, {1}/{\sqrt{2}}, 0, \dots, 0) \mbox{ and } \nu =({1}/{\sqrt{2}}, -{1}/{\sqrt{2}}, 0, \dots, 0).
$$
Then 
$$
V_{n-1}(O_n | \xi^\perp) =V_{n-1}(O_n | \nu^\perp) =\frac{1}{\sqrt{2}}V_{n-1}(O_{n-1})=\frac{1}{\sqrt{2}}\frac {2^{n-1}}{(n-1)!}.
$$
Also notice that $V_{n-2}(O_n |  (\xi, \nu)^\perp)=V_{n-2}(O_{n-2})=2^{n-2}/(n-2)!$. Thus inequality (\ref{proj}) is false when $D=O_n$.

We next show that (\ref{proj}) is false if $D$ is a rotation invariant body, which gives another evidence of Conjecture \ref{Cj:BE}.
\begin{Prop}
 Consider a concave function $f:[a,b] \to \R^+$ and define convex body 
$
K_f=\{x \in \R^n\ |\ \sum_{i=2}^n x_i^2 \le f^2(x_1)\}.
$
Choose $\xi =(0, 1,0, \dots, 0) \mbox{ and } \nu =(0,0, 1,0, \dots, 0)$.  Then 
$$
\frac{n}{n-1}V_{n-2}(K_f |  (\xi, \nu)^\perp) V_n(K_f)>V_{n-1}(K_f| \xi^\perp) V_{n-1}(K_f | \nu^\perp).
$$
\end{Prop}
\begin{pf} Using the Fubini theorem, we get
$$
V_n(K_f)=\int_{a}^b f(t)^{n-1} V_{n-1}(B^{n-1}) dt=\kappa_{n-1}\int_{a}^b f(t)^{n-1} dt,
$$
where $B^{n} =\{x \in \R^n\ |\ \sum x_i^2 \le 1\} $ is the Euclidean unit ball and $\kappa_n=V_n(B^n)$. Moreover,
$$
V_{n-1}(K_f | \xi^\perp )=V_{n-1}(K_f | \nu^\perp )=\kappa_{n-2}\int_{a}^b f(t)^{n-2} dt,
$$
and
$$
V_{n-2}(K_f | (\nu, \xi)^\perp )=\kappa_{n-3}\int_{a}^b f(t)^{n-3}dt.
$$
Thus to prove the proposition we need to show that 
$$
\left(\kappa_{n-2}\int_{a}^b f(t)^{n-2}  dt \right)^2<
\frac{n}{n-1} \left( \kappa_{n-3}\int_{a}^b f(t)^{n-3} dt \right) \left(\kappa_{n-1}\int_{a}^b f(t)^{n-1} ) dt \right).\nonumber
$$
First, by H\"older's inequality, we get
$$
\left(\int_{a}^b f(t)^{n-2} dt \right)^2 \le \left( \int_{a}^b f(t)^{n-3} dt \right) \left(\int_{a}^b f(t)^{n-1} dt \right).
$$
Next, we will show that
$$
\kappa_{n-2}^2 < \frac{n}{n-1} \kappa_{n-3}\kappa_{n-1}.
$$
Indeed,  using $\kappa_k=\pi^{k/2}/\Gamma(k/2+1)$ (see \cite{MS}, \cite{Sch1}) we get that the above inequality is equivalent to
$$
\Gamma^{2}(n/2)> \frac{n-1}{n} \Gamma(n/2 -1/2)\Gamma( n/2+1/2),
$$
which follows immediately from the classical inequality for the Gamma function:
$$
\sqrt{x+\frac{1}{2}} > \frac{\Gamma(x+1)}{\Gamma(x+1/2)} .
$$
\end{pf}

\section{Isomorphic Version}\label{S:iso}

As we have already seen, \re{BE} may not hold for a general convex body $D$. However, we can relax the inequality by introducing a constant  $c_{n,r} >0$ depending on the dimension $n$ and the
number $r$ only:
\begin{equation}\label{e:isoBE2}
V(K_1,\dots, K_r,D^{n-r})V_n(D)^{r-1}\leq c_{n,r} \prod_{i=1}^r V(K_i,D^{n-1}).
\end{equation}
It is an interesting question to find the minimal such constant $c_{n,r}$. In this section we give 
an upper bound for $c_{n,r}$ in the general and in the symmetric cases (see \rt{general}). 
In addition, we prove that $c_{n,r}= r^{r}/r!$ in the special case when $K_1,\dots, K_r\subset\R^n$ are zonoids and $D\subset\R^n$ is an arbitrary convex body
(see \rt{zonoids}). 

We start with a rough estimate for the constant $c_{n,r}$.

\begin{Prop}\label{P:general} There exists a constant $c_{n,r}<(n+2)^{n-r} n^{n(r-1)}$ such that
$$
V(K_1,\dots, K_r,D^{n-r})V_n(D)^{r-1}\leq c_{n,r} \prod_{i=1}^r V(K_i,D^{n-1}).
$$
is true for all convex bodies $K_1,\dots, K_r,$ and $D$ in~$\R^n$.
\end{Prop}
\begin{pf} It is enough to consider the case when $D$ is $n$-dimensional.
We use a standard idea of approximating $D$ by a simplex of maximal volume contained in $D$ (see  \cite{Gru1},  \cite{L1}, \cite{L2}  and \cite{Sch2}).
Let $\D$ be a simplex of maximal volume inscribed in~$D$. Using the invariance of (\ref{e:isoBE2}) under affine transformation we may assume that $\D$ is a regular simplex with the barycenter at 
the origin.  For every vertex  $v \in \D$ choose a hyperplane $H$ such that  $v \in H$ and  $H$ is parallel to the facet of $\D$ opposite to $v$. Then $H$ supports $D$, otherwise we would get a contradiction with the maximality of $\D$. This shows that $D \subset -n \D$. (Note that if  $D$ is a symmetric convex body then we immediately get $D \subset n\D$, which would slightly improve our estimate on the  constant $c_{n,r}$.) Also by a result of Lassak \cite{L1} $ D\subseteq(n+2)\D.$  Therefore we obtain
$$
\D \subset D\subseteq\left((n+2)\D\right) \cap \left(-n \D\right).
$$
From Theorem \ref{T:ineq} we get that
$$
V(K_1,\dots, K_r,\D^{n-r})V_n(\D)^{r-1}\leq \prod_{i=1}^r V(K_i,\D^{n-1}),
$$
and thus
$$
(n+2)^{-(n-r)}V(K_1,\dots, K_r,D^{n-r})  n^{n(1-r)}V_n(D)^{r-1}\leq \prod_{i=1}^r V(K_i, D^{n-1}).
$$

\end{pf}

One can improve the above estimate   for  $c_{n,r}$ by
 computing the volume of $\left((n+2)\D\right) \cap \left(-n \D\right)$. This would show 
$c_{n,r}<(c n)^{r(n-1)}$ for some absolute constant $c<1$, which we still feel is not the optional bound.
Instead, we are going to first prove a better estimate in the case of $K_i$ being zonoids and 
then apply John's Theorem (\cite{J}, see also \cite{MS}).

First recall the inequality for projections in \cite[Lemma 4.1]{GHP}  (see also \cite{FGM}):
$$
 \frac{n}{2(n-1)}V_{n-2}(D |  (\xi, \nu)^\perp) V_n(D) \le V_{n-1}(D | \xi^\perp) V_{n-1}(D | \nu^\perp)
$$
for any convex body $D\subset\R^n$. Thus using (5.3.23)  of \cite[page 294]{Sch1}  we see that
\begin{equation}\label{e:proj}
V(K_1,K_2,D^{n-2})V_n(D)\leq 2 V(K_1,D^{n-1})V(K_2,D^{n-1})
\end{equation}
is true for any convex body $D$ and two orthogonal segments $K_1,K_2$. Second, the fact that the above inequality is invariant under linear transformations implies that the orthogonality assumption is not necessary. 
Finally, by the additivity of the mixed volume with respect to the Minkowski sum and continuity with respect
to the Hausdorff metric, \re{proj} generalizes to arbitrary zonoids $K_1, K_2$. We have obtained the following lemma.

\begin{Lemma} For any zonoids $K_1, K_2 \subset \R^n$ and any convex  body $D\subset\R^n$ we have
$$
V(K_1,K_2,D^{n-2})V_n(D) \leq 2V(K_1,D^{n-1})V(K_2,D^{n-1}).
$$
\end{Lemma}

It is well known that any symmetric convex  body $D$ in $\R^2$ is a zonoid (see \cite[Theorem 3.2.11]{Sch1}). 
This implies the following.

\begin{Cor}\label{C:symm-plane} Suppose $K_1$, $K_2$ and $D$ are convex bodies in $\R^2$  and  
$K_1,K_2$ are symmetric. Then
$$
V(K_1,K_2)V_2(D) \leq 2 V(K_1,D)V(K_2,D).
$$
\end{Cor}
The next fact is an analog of \cite[Lemma 4.1]{GHP}.

\begin{Lemma} \label{lm:proj} Consider a convex body $D\subset\R^n$. Let  $P_i(D)$ be the orthogonal projection of $D$ onto $e_i^\perp$ and $P_{[r]} (D)$  be the orthogonal projection of $D$ onto
$(e_1, e_2, \dots, e_r)^\perp$. Then
$$
\left(\frac{n}{r}\right)^r \binom{n}{r}^{-1} V_{n-r}(P_{[r]} (D)) V_n(D)^{r-1} \le \prod\limits_{i=1}^r V_{n-1}(P_i(D)).
$$
\end{Lemma}
\begin{pf}   Consider $x \in \R^n$ and write $x=(x_1,\dots, x_r, y)$ where $y \in \R^{n-r}$. For every $y \in P_{[r]} (D)$ define
$$
D_i(y) =\{(x_1,\dots, x_{i-1}, x_{i+1}, \dots, x_r)\ |\  (x_1,\dots, x_{i-1}, 0, x_{i+1}, \dots, x_r, y) \in P_i(D)\}
$$
and
$$
D_{[r]} (y)=\{(x_1, \dots, x_r) \in \R^r\ |\  (x_1, \dots, x_r, y) \in D\}.
$$
Note that $D_i(y)$ is the orthogonal projection of  $D_{[r]} (y)$ onto the coordinate plane $e_i^\perp$ in $\R^r$ and thus by the Loomis-Whitney inequality  \cite{LW}
$$
V_{r}(D_{[r]} (y))^{r-1} \le \prod\limits_{i=1}^r V_{r-1}(D_i(y)).
$$
Next
$$
 \prod\limits_{i=1}^r V_{n-1}(P_i(D))  =  \prod\limits_{i=1}^r\, \int\limits_{P_{[r]} (D)} V_{r-1}(D_i(y)) dy\ge$$
 $$
\left[\,\int\limits_{P_{[r]} (D)} \prod\limits_{i=1}^r  V^{\frac{1}{r}}_{r-1}(D_i(y)) dy \right]^r\ge \left[\,\int\limits_{P_{[r]} (D)}   V_{r}(D_{[r]} (y))^{\frac{r-1}{r}} dy \right]^r.
$$
By the Brunn-Minkowski inequality the function $\phi(y)=V_{r}(D_{[r]} (y))^{\frac{1}{r}}$ is concave on  $P_{[r]} (D)$. We apply Berwald's Lemma (\cite{B}, see also \cite[Lemma 4.2]{GHP}) with 
$$
\phi(y)=V_{r}(D_{[r]} (y))^{\frac{1}{r}},
$$
$$
A = P_{[r]} (D);\ \  p=r-1,\ \  q=r,\  \mbox{ and }\  k=n-r.
$$
We get
$$
\left[ \binom{n-1}{n-r}  \frac{1}{ V_{n-r} (P_{[r]} (D))} \int\limits_{P_{[r]} (D)} |\phi(y)|^{r-1} \right]^{\frac{1}{r-1}}
\ge
\left[ \binom{n}{n-r}  \frac{1}{ V_{n-r} (P_{[r]} (D))} \int\limits_{P_{[r]} (D)} |\phi(y)|^{r} \right]^{\frac{1}{r}}$$
which is equivalent to 
$$
\left[  \int\limits_{P_{[r]} (D)} |\phi(y)|^{r-1} \right]^{r}
\ge
V_{n-r} (P_{[r]} (D))\frac{\binom{n}{n-r}^{r-1}}{\binom{n-1}{n-r}^r}\left[  \int\limits_{P_{[r]} (D)} |\phi(y)|^{r} \right]^{r-1},
$$ i.e.
$$
\left[  \int\limits_{P_{[r]} (D)} V_{r}(D_{[r]} (y))^{\frac{r-1}{r}} \right]^{r}
\ge
V_{n-r} (P_{[r]} (D))\frac{\binom{n}{n-r}^{r-1}}{\binom{n-1}{n-r}^r}\left[  \int\limits_{P_{[r]} (D)} V_{r}(D_{[r]} (y))\right]^{r-1},$$
and the statement of the lemma follows.
\end{pf}
\begin{Rem}\label{rem:zon} We note that the estimate in Lemma \ref{lm:proj} is the best possible. Indeed, let  $O^k$ be the octahedron as in \rs{proj} and $Q^k=\{x\in \R^k\ |\ |x_i| \le 1\}$ a cube in $\R^k$. 
We define the convex body $D$ as the convex hull of 
the union $\left(Q^r\times\{0\}\right)\cup\left(\{0\}\times O^{n-r}\right)\subset \R^r\times\R^{n-r}=\R^n$. Then
$$
V_n(D)=2^n \frac{r !}{n!}.
$$
Furthermore, $P_{[r]} (D) =O^{n-r}$  and if $1\leq i\leq r$ then  $P_i(D)$ is the convex hull of the union of
an $(r-1)$-dimensional cube and the octahedron $\{0\}\times O^{n-r}$.
Thus
$$
V_{n-1}(P_i(D))=2^{n-1} \frac{(r-1)!}{(n-1)!}\quad \mbox{and} \quad V_{n-r}(P_{[r]} (D))=\frac{2^{n-r}}{(n-r)!}.
$$
Therefore $D$ gives equality in the statement of Lemma \ref{lm:proj}.
\end{Rem}
Now let us restate Lemma \ref{lm:proj} in the language of mixed volumes. Let 
$K_1,\dots, K_r$ be pairwise orthogonal unit segments. Then 
again by (5.3.23)  in \cite[page 294]{Sch1},
$$
V(K_i, D^{(n-1)})=\frac{1}{n}V_{n-1}(D | K_i^\perp).
$$
Moreover, if $U=\sum_{i=1}^r K_i$ then
$$
\binom{n}{r}V(D^{n-r}, U^r)=V_r(U) V_{n-r}(D| (K_1, \dots, K_r)^\perp).
$$
Thus
$$
V_{n-r}(D| (K_1, \dots, K_r)^\perp) =\binom{n}{r}r!V(K_1, \dots, K_r, D^{n-r}),
$$
and  Lemma \ref{lm:proj} implies
$$
\frac{r!}{r^{r}} V(K_1, \dots, K_r, D^{n-r})V_n(D)^{r-1} \le \prod\limits_{i=1}^r V(K_i, D^{n-1}).
$$
Since the above formula is invariant under a linear transformation the orthogonally 
assumption on $K_1, \dots, K_r$ is not necessary. Again, by the additivity and continuity of the mixed volume,
and taking into account the example in Remark \ref{rem:zon} we obtain the following result.

\begin{Th}\label{T:zonoids} Suppose $D$ is a convex body in $\R^n$  and   $K_1,\dots K_r$ are zonoids  then
$$
 V(K_1, \dots, K_r, D^{n-r})V_n(D)^{r-1} \le \frac{r^{r}}{r!} \prod\limits_{i=1}^r V(K_i, D^{n-1}),
$$
and the inequality is sharp.
\end{Th}

Now we may apply John's theorem to improve the bound for $c_{n,r}$ 
in \rp{general}. Indeed, for any convex body $K_i$ there exists an ellipsoid $E_i$ such that
$E_i \subset K_i \subset n E_i$. Moreover,  $E_i \subset K_i \subset \sqrt{n} E_i$ when $K_i$ is symmetric. It is also well known that an ellipsoid is a particular case of a zonoid. This implies the following theorem.

\begin{Th}\label{T:general} There exists a constant $c_{n,r}\le n^rr^{r}/ r!$ such that
$$
V(K_1,\dots, K_r,D^{n-r})V_n(D)^{r-1}\leq c_{n,r} \prod_{i=1}^r V(K_i,D^{n-1})
$$
holds for all convex bodies $K_1,\dots, K_r$ and $D$ in $\R^n$. Moreover $c_{n,r}\le {n}^{r/2}r^{r}/ r!$ when $K_1,\dots, K_r$ are symmetric with respect to the origin.
\end{Th}

\section{Isomorphic Bezout inequality in dimension 2}\label{S:dim2}

Applying \rt{general} in the case of $n=r=2$ and $K_1,K_2$ being symmetric we obtain $c_{2,2}\leq 2$ (see also \rc{symm-plane}). It turns out  that $c_{2,2}= 2$ always, even in the non-symmetric case.

\begin{Th}\label{T:plane} Let $K$, $L$  and $D$ be any convex bodies in $\R^2$. Then
\begin{equation}\label{e:plane}
V(K,L)V_2(D)\leq 2V(K,D)V(L,D).
\end{equation}
Moreover, the equality is attained when $K$ and $L$ are line segments and $D=\lambda_1K+\lambda_2L$ for
some $\lambda_i\geq 0$.
\end{Th}
\begin{pf} The second statement follows immediately since $V_2(D)=2\lambda_1\lambda_2V(K,L)$
and $V(K,D)=\lambda_1V(K,L)$,  $V(L,D)=\lambda_2V(K,L)$, by the additivity of the mixed volume.

To show the inequality \re{plane}, first note that it is enough to consider the case of convex polygons $K,L$
as they approximate convex 2-dimensional bodies.
Recall that any convex polygon is the Minkowski
sum of triangles and line segments (see \cite[p. 177]{YB}). Let 
$$K=\sum_{i=1}^k K_i\quad L=\sum_{i=1}^l L_i$$
be the corresponding Minkowski decompositions into triangles and line segments. Then, by additivity, \re{plane}
can be written as 
$$\sum_{i=1}^k \sum_{i=1}^lV(K_i,L_j)V_2(D)\leq\sum_{i=1}^k \sum_{i=1}^l 2V(K_i,D)V(L_j,D).$$
This reduces the proof to the case of two triangles $K$ and $L$. (The case when one or both
of $K,L$ is a line segment follows by contracting one of the sides of a triangle to a point.)
Furthermore, \re{plane} is invariant under translations and dilations of each $K$, $L$, and $D$.
Thus, after a possible dilation and translation of $K$ and $L$ we may assume that $D$ is inscribed in both $K$
and $L$, in which case $V(K,D)=V_2(K)$ and $V(L,D)=V_2(L)$ by \rc{cor}. Therefore, we may 
assume that  $D=K\cap L$. The statement now follows from \rt{triangles} below.
\end{pf}

In the statement and the proof of the next theorem we use $|K|$ to denote the area (2-dimensional volume) of 
a triangle $K$ and $|\overline{AB}|$ the length of a line segment $\overline{AB}$.

\begin{Th}\label{T:triangles}
Let $K,L$ be two triangles in the plane whose intersection is inscribed in each of them. Then 
\begin{equation}\label{e:triangles}
V(K,L)|K\cap L|\leq 2|K||L|.
\end{equation}
\end{Th}
\begin{pf}
Let $c>0$ be the smallest number such that $K$ is contained in $L'=cL+v$ for some $v\in \R^2$.
We distinguish two cases: either the normals to the sides of $K$ and $L$ are distinct and
alternate (as points on $\mS^1$) or not. Notice that in the latter case $K$ and $L'$ have a common vertex,
and in the former the three vertices of $K$ lie on three different sides of $L'$.
Note that \re{triangles} is invariant under affine transformations applied to both $K$ and $L$.
Therefore we may assume that $L'$ equals
the standard $2$-simplex $\D$. 

First, suppose $K$ and $L'$ share a vertex, which we may assume to be $(1,0)$.
Then $K=\triangle ABC$ with $C=(1, 0)$, and  
$L=\triangle A'B'C'$ is a translate of $\mu\D$ 
where $\mu=1/c>0$ (see \rf{triangles1}). 
\begin{figure}[h]
\includegraphics[scale=1.9]{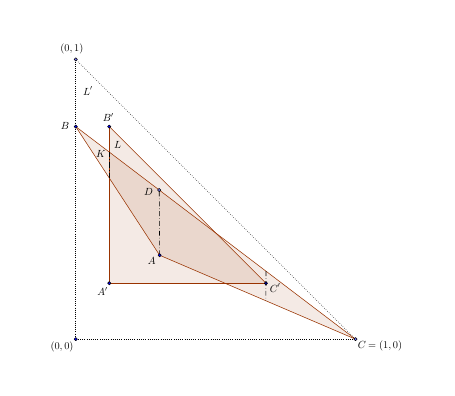}
\caption{Case 1.}
\label{F:triangles1}
\end{figure}

Note that $K$ is inscribed in $\D$, so by \rc{cor},
$$V(K,L)=V(K,\mu\D)=\mu V(K,\D)=\mu|\D|=\mu/2.$$ 
Since $|L|=\mu^2/2$  the inequality \re{triangles} becomes
$$|K\cap L|\leq 2\mu|K|.$$

Consider the vertical segment $\overline{AD}$ subdividing $K$ into two triangles as in \rf{triangles1}.
Note that its length $|\overline{AD}|$ equals $2|K|$ divided by the horizontal width of $K$, which
by our choice of coordinates equals 1, i.e. $|\overline{AD}|=2|K|$. 
Next, notice that $K\cap L$ is contained in the vertical strip of width $\mu$ and its 
longest vertical section is at most $|\overline{AD}|$. Therefore,
$$|K\cap L|\leq \mu|\overline{AD}|= 2\mu|K|,$$
and this case is proved. 

For the other case, $K$ has vertices $A=(\alpha,0)$, $B=(0,\beta)$,  $C=(\gamma, 1-\gamma)$ for
some $0<\alpha,\beta,\gamma < 1$, and
$L=\triangle A'B'C'$ is a translate of $\mu\D$ (see \rf{triangles}).
\begin{figure}[h]
\includegraphics[scale=1.3]{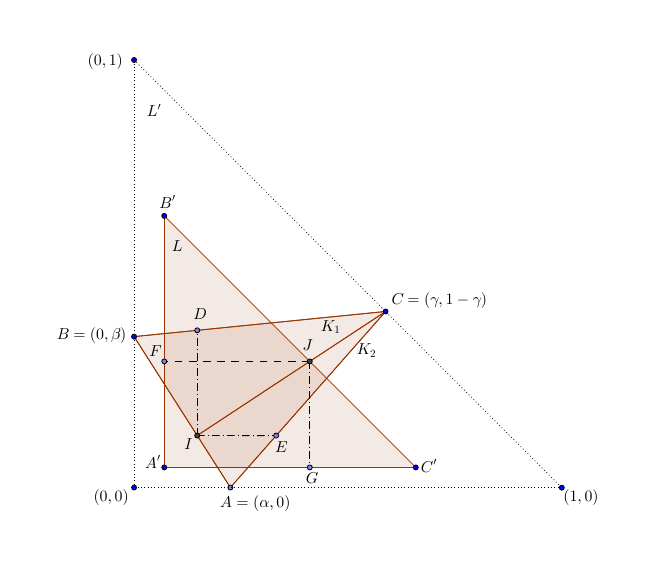}
\caption{Case 2.}
\label{F:triangles}
\end{figure}
As in the previous case, $K$ is inscribed in $\D$, so we need to show $$|K\cap L|\leq 2\mu|K|.$$
We subdivide $K$ into two triangles $K_1=\triangle BIC$ and
$K_2=\triangle AIC$ where $I$ is a point on the side $\overline{AB}$ chosen in such a way that
$|K_1|/|K|=\gamma$.
In coordinates
$$I=(t\alpha, (1-\gamma)\beta).$$ 

Next, let $\overline{ID}$ and $\overline{IE}$ be the vertical and the horizontal segments
subdividing the triangles $K_1$ and $K_2$, respectively. Using the same argument as above we get 
$$\frac{2|K_1|}{\gamma}=|\overline{ID}|,\quad \text{and}\quad \frac{2|K_2|}{1-\gamma}=|\overline{IE}|.$$
Furthermore, let $J$ be the intersection of $\overline{CI}$ and $\overline{B'C'}$, and let $\overline{JF}$ and 
$\overline{JG}$ be
the horizontal and the vertical segments as in \rf{triangles}. Then, similar to
the first case, we have 
$$|K_1\cap L|\leq |\overline{ID}||\overline{JF}|\quad\text{and}\quad |K_2\cap L|\leq |\overline{IE}||\overline{JG}|.$$
We denote $\mu_1=|\overline{JF}|$, and so $\mu-\mu_1=|\overline{JG}|$.
Therefore we can write
$$|K\cap L|\leq \frac{2|K_1|}{\gamma}\mu_1+\frac{2|K_2|}{1-\gamma}(\mu-\mu_1).$$
It remains to show that 
$$\frac{|K_1|}{\gamma}\mu_1+\frac{|K_2|}{1-\gamma}(\mu-\mu_1)\leq \mu|K|,$$
which is immediate once we note that $|K_1|/|K|=\gamma$ and $|K_2|/|K|=1-\gamma$.
\end{pf}

As a corollary we obtain the following interesting result in elementary plane geometry.

\begin{Cor} Let $H=AA'BB'CC'$ be a convex hexagon and $K=\triangle ABC$, $L=\triangle A'B'C'$.
Then 
$$|H||K\cap L|\leq 2|K||L|.$$
\end{Cor}

\begin{pf} This follows from \rt{triangles} and the fact that $V(K,L)=|H|$. Indeed,
$|H|=V(H,L)$, since  $L$ is inscribed in $H$ (see \rc{cor}).
On the other hand, 
$$V(K,L)=\frac{1}{2}\sum_{u}h_K(u)|L^u|,$$
where the sum is over the unit normals to $L$. Clearly, $h_H(u)=h_K(u)$ for all such $u$ and
so $V(K,L)=V(H,L)$.

\end{pf}

\end{document}